\documentclass[11pt,twoside]{article}
\usepackage{datetime}\usdate
\usepackage{graphicx}
\usepackage{artmods,macros}
\usepackage{amsmath, amssymb}
\usepackage{multirow}
\usepackage[
backend=bibtex,
style=numeric,
sorting=nyt,
giveninits=true
]{biblatex}
\usepackage{algorithm}
\usepackage{algpseudocode}

\date{}

\title{Finite Newton-Step Generation via Subspace Quasi-Newton Methods for Quadratic Problems\thanks{\footWASP}}
\author{Aban ANSARI-\"{O}NNESTAM\thanks{\footAO} \and Anders FORSGREN\thanks{\footAF}}

\def\footAF{Division of Numerical Analysis, Optimization and Systems
  Theory, Department of Mathematics, KTH Royal Institute of Technology, SE-100 44
            Stockholm, Sweden ({\tt andersf@kth.se}).}

 \def\footAO{Division of Applied Mathematics, Department of
            Mathematics, Link\"{o}ping University, SE-581 83
            Link\"{o}ping, Sweden ({\tt aban.ansari-onnestam@liu.se}). }
\def\footWASP{This work was partially supported by the Wallenberg AI, Autonomous Systems and
Software Program (WASP) funded by the Knut and Alice Wallenberg Foundation.}

   \input macros

\def\0{^{(0)}}
\def\1{^{(1)}}

\def\st{\mbox{\rm subject to}}

\def\rn{\mathbb{R}^n}
\def\K{\mathcal{K}}

\usepackage{color}

\defbibheading{bibintoc}[\refname]{%
  \section*{#1}
  \markboth{#1}{#1}
}

\begin{document}

\maketitle\thispagestyle{empty}

\begin{abstract}
We study the curvature information sufficient for finite Newton-step generation by quasi-Newton methods on strictly convex quadratic problems.
In exact arithmetic, the memoryless BFGS quasi-Newton method terminates
finitely on such a problem when exact line search
is used, although its quasi-Newton approximation incorporates
curvature information associated with only one direction. We show that
exact line search can be replaced by curvature information associated
with one additional direction. Information associated with at most two
directions suffices to generate the Newton step after finitely many
iterations, independently of the preceding step lengths. More precisely, knowledge of the action of the Hessian on a suitably chosen subspace of dimension at most two is sufficient to generate a step that decomposes into a subspace Newton step and a component parallel to the next conjugate direction. Applying this
construction recursively generates the full Newton step after finitely
many iterations, after which a unit step reaches the minimizer. The
required Hessian actions can also be recovered from gradient
differences, yielding a first-order formulation. The resulting algorithm is
intended as a constructive device for identifying sufficient curvature
information, rather than as a practical alternative to the method of
conjugate gradients.
\end{abstract}

\begin{keywords}
quasi-Newton method, unconstrained quadratic problem, finite termination, Krylov subspace method.
\end{keywords}

\section{Introduction}
We consider the strictly convex quadratic problem
\begin{equation}\label{qp}
  \begin{array}{ll}
    \minimize{x \in \rn} & \frac{1}{2}x^{T}Hx + c^{T}x 
\end{array}
\tag{QP}
\end{equation}
where $H=H^{T}\succ0$. Its unique minimizer $x^*$ is characterized by the linear system $Hx^*+c=0$.

Given an initial point $x_0$, a method using mutually $H$-conjugate search directions and exact line search reaches the minimizer of \eqref{qp} in at most $n$ iterations. At iteration $k$, the next iterate is
\[
 x_{k+1}=x_k+\alpha_kp_k,
\]
where $\alpha_k$ is obtained by exact line search along $p_k$. In a
quasi-Newton method, $p_k$ is determined from $B_kp_k=-g_k$, where
$g_k=g(x_k)$ represents the gradient at $x_k$ and $B_k=B_k^{T}\succ0$ approximates $H$.

This paper investigates what curvature information is sufficient to generate the Newton step finitely when exact line search is not used. Under
exact line search, the memoryless BFGS quasi-Newton method achieves
finite termination while incorporating curvature information
associated with only the most recent direction. We ask what curvature
information is sufficient when exact line search is removed and the
preceding step lengths are arbitrary. Our central observation is that
the quasi-Newton approximation need not agree with $H$ on the entire
Krylov subspace. Agreement on a subspace generated by at most two
appropriately chosen vectors is sufficient. This yields an iterative
construction that allows arbitrary preceding step lengths and
generates the Newton step after finitely many iterations. Once that
step has been generated, a unit step reaches the minimizer.

Section~\ref{sec:background} reviews Krylov-subspace and quasi-Newton methods. Section~\ref{sec:subspace} introduces subspace Newton steps and describes their evolution as the Krylov subspace expands. Section~\ref{sec:qnsubspace} develops the corresponding subspace quasi-Newton construction and shows that information associated with at most two vectors is sufficient. Section~\ref{sec:qn} gives a recursive realization of the construction for arbitrary preceding step lengths.

\section{Background}\label{sec:background}

Using $H$-conjugate search directions and exact line searches, one can solve \eqref{qp} in at most $n$ iterations.
Given a point $x \in \rn$ and a search direction $p$, exact line search can be used to calculate the optimal step length
\[ \alpha = -\frac{g(x)^{T}p}{p^{T}Hp}. \]
We refer to this as the \textit{Newton scaling of $p$ with respect to $x$}.
Fix $x_0 \in \rn$, and suppose that $\{p_0, \dots, p_{n-1}\}$ forms a set of conjugate vectors with respect to $H$, that is, $p_i^{T}Hp_j=0$ for $i\neq j$.
If exact line search is used to calculate the step size $\alpha_k$, then $x_{k+1} = x_k + \alpha_kp_k$ is the solution to the constrained optimization problem
\[
  \begin{array}{ll}
   \minimize{x \in \rn} & \frac{1}{2}x^{T}Hx + c^{T}x \\[4pt]
   \st & x \in x_0 + \operatorname{span} \{p_0, \dots, p_k\}.
\end{array}
\]
We focus on Krylov-subspace methods, in which the search directions span a Krylov subspace.

\begin{definition}
Given a vector $b$ and a matrix $A$, the $k$-th Krylov subspace generated by $b$ and $A$ is
\begin{equation*}
\K_k(b,A) = \operatorname{span} \{b, Ab, \dots, A^{k-1}b\}.
\end{equation*}
We use the term \emph{Krylov-subspace method} for a
method that successively generates the minimizers over the expanding
affine spaces $x_0 + \K_k(g_0,H)$.
\end{definition} 

\begin{definition}
  Given $x_0 \in \rn$, let $r$ be the smallest integer such that $\K_r(g_0,H) = \K_{r+1}(g_0,H)$. We adopt the conventions $\K_0(g_0,H)=\{0\}$, $\widehat{x}_0=x_0$, and $\widehat{g}_0=g_0$.
  For $k \leq r$, denote the minimizer of \eqref{qp} over the affine space $x_0 + \K_k(g_0,H)$ as $\hat x_k$. The gradient at the minimizer is denoted by $\hat g_k$ and the search direction generated by such a Krylov-subspace method is parallel to the difference $\hat q_{k-1} = \hat x_k - \hat x_{k-1}$.
\end{definition}

The search directions $\hat q_0, \dots, \hat q_{k-1}$ are mutually conjugate with respect to $H$ and form a basis for $\K_k(g_0,H)$. Additionally, the gradients $\hat g_0, \dots, \hat g_{k-1}$ form an orthogonal basis for the same subspace.
For $r$ as defined above, the minimizer $\hat x_r$ is exactly the global minimizer $x^*$ of \eqref{qp}. For a more comprehensive overview of Krylov-subspace methods, see, e.g., \cite[Chapter 5]{nocedal}, \cite[Chapter 6]{Saa95}.

\subsection{Generating conjugate search directions by exact line search}\label{subsec:conjugate}

In order to solve \eqref{qp} using a Krylov-subspace method, search directions that are parallel to $\hat q_k$ must be generated. Parallel vectors will be denoted by $p \parallel p'$.
When exact line search is used on these parallel directions, the iterates are identical, and as such after $r$ iterations the algorithm will terminate.

Conjugate directions can be generated by the conjugate Gram--Schmidt process.
The method of conjugate gradients (CG) generates these search directions iteratively by setting
\begin{gather*}
p_0 = -g_0, \quad
p_k = -g_k + \frac{g_k^{T}Hp_{k-1}}{p_{k-1}^{T}Hp_{k-1}}p_{k-1}
= -g_k + \frac{g_k^{T}g_{k}}{g_{k-1}^{T}g_{k-1}}p_{k-1}.
\end{gather*}
In the classical CG method, exact line search ensures that each
gradient is orthogonal to the preceding search directions and that the
generated directions are mutually $H$-conjugate, yielding finite
termination in exact arithmetic. For a detailed description of the
method of conjugate gradients, see, e.g., Saad~\cite[Chapter 6.7]{Saa95}.

\subsection{Quasi-Newton methods}\label{subsec:qnmethods}

For a quadratic objective, the Broyden--Fletcher--Goldfarb--Shanno (BFGS) method also generates $H$-conjugate directions when exact line search is used.
Given an initial Hessian approximation $B_0 \succ 0$, at each iteration a search direction is obtained by solving the system $B_k p_k = -g_k$ and the matrix is updated by
\[ B_{k+1} = B_k - \frac{1}{p_k^{T}B_kp_k}B_k p_k p_k^{T} B_k + \frac{1}{p_k^{T}Hp_k} H p_k p_k^{T} H.\]

This can be extended to a memoryless variant where the Hessian approximation is built using the identity matrix instead of the previous approximation $B_k$, so that
\[ B_{k+1} = I - \frac{1}{p_k^{T}p_k}p_k p_k^{T} + \frac{1}{p_k^{T}Hp_k} H p_k p_k^{T} H.\]

Under exact line search, CG, BFGS, and memoryless BFGS generate directions $p_k \parallel \hat q_k$ and are therefore Krylov-subspace methods \cite{huang}.

The first quasi-Newton method was proposed by
Davidon~\cite{davidon} and subsequently developed by Fletcher
and Powell~\cite{fletcher}. For an overview on quasi-Newton methods,
see, e.g., \cite[Chapter 4]{nocedal}, \cite[Chapter 4.5]{GMW19} and
\cite[Chapter 3.5]{Fle81}.  There is a vast literature on quasi-Newton
methods, which is beyond the scope of the current manuscript. We focus
on the BFGS method and the memoryless BFGS method, which are of
relevance for our results.

Under suitable conditions and with exact line search, quasi-Newton
methods generate directions parallel to the corresponding Krylov-space
directions and terminate in exactly $r$ iterations
\cite{forsgren2015connection,forsgren2018exact}.

Kolda, O'Leary and Nazareth \cite{kolda} show that finite termination
remains possible even if the quasi-Newton matrix is not updated at
every iteration, provided that exact line search is
used. Powell~\cite{powell1, powell2} similarly proves finite
termination using unit steps as long as the update of the quasi-Newton
matrix $B_k$ is skipped in every other iteration.

\section{Subspace Newton steps}\label{sec:subspace}
Suppose for some initial point $x_0$ we consider the Krylov subspace $\K_{k}(g_0,H)$.
Then the Krylov-subspace minimizer $\hat x_k$ will belong to the affine space $x_0 + \K_{k}(g_0, H)$.
For any point $x \in x_0 + \K_{k}(g_0,H)$ there exists a unique \textit{subspace Newton step} $p_{k-1}^{N}(x) \in \K_k(g_0,H)$ such that $x + p_{k-1}^{N}(x) = \hat x_k$.
Let $S=[s_0\ \cdots\ s_{k-1}]$, where the columns of $S$ form a basis for $\K_{k}(g_0,H)$. If $g(x)$ is the gradient at the given point $x$, then the subspace Newton step $p_{k-1}^{N}(x)$ can be expressed as
\[
    p_{k-1}^{N}(x) = S\beta,
\]
for $\beta$ given by
\[
S^{T}HS\beta = -S^{T}g(x).
\]

This follows from the fact that $\hat x_k$ is the unique point that minimizes \eqref{qp} over the affine space $x_0 + \K_k(g_0,H)$, or equivalently that the gradient at the minimizer $\hat g_k$ is orthogonal to $\K_{k}(g_0,H)$.
The subspace Newton step is fundamental to our analysis because the relative change of the directional derivatives
  of the gradient along the subspace Newton step with respect to all
  conjugate directions can be characterized by a single step length
  $\alpha$, as stated in the following lemma.

\begin{lemma}\label{newton:step}
  For $k \leq r$, let $q_0, \dots, q_{k-1}$ be a basis for $\K_k(g_0,H)$ where $q_i \parallel \hat q_i$. For any $x \in x_0 + \K_k(g_0,H)$, the subspace Newton step $p_{k-1}^{N}(x)$ to $\hat x_k$ is
  \[
    p_{k-1}^N(x) = \sum_{i=0}^{k-1}-\frac{g(x)^{T}q_i}{q_i^{T}Hq_i}q_i.
  \]
  Additionally, for any $\alpha \in \mathbb{R}$
  \[
    g(x+\alpha p_{k-1}^N(x))^{T} q_j =(1-\alpha) g(x)^{T} q_j, \quad j=0,\dots,k-1.
  \]
\end{lemma}

\begin{proof}
Since $q_i \parallel \hat q_i, i = 0, \dots, k-1$, the basis $q_0, \dots, q_{k-1}$ is mutually conjugate with respect to $H$.
We may therefore represent the subspace Newton step as $p_{k-1}^{N}(x) = \sum_{i=0}^{k-1}\beta_i q_i$ where $\beta_i, i = 0, \dots, k-1$ are the solutions to the system
\[
  \mtx{ccccc}{ q_{k-1}^{T} H q_{k-1} & 0 & 0 & 0 \\
                      0  & q_{k-2}^{T} H q_{k-2} & 0 & 0 \\
                      0 & 0 &  \ddots & 0 \\
                      0 & 0 & 0 & q_0^{T} H q_0 }
  \mtx{c}{\beta_{k-1} \\ \beta_{k-2} \\ \vdots \\ \beta_0 }
  =
  -   \mtx{c}{q_{k-1}^{T} g(x) \\ q_{k-2}^{T} g(x) \\ \vdots \\
    q_0^{T} g(x) }.
\]
Therefore it is straightforward to show that, \[\beta_i = -\frac{g(x)^{T}q_i}{q_i^{T}Hq_i}, \qquad i = 0, \dots, k-1\]
and the subspace Newton step is
\[p_{k-1}^N(x) = \sum_{i=0}^{k-1}-\frac{g(x)^{T}q_i}{q_i^{T}Hq_i}q_i.
\]
  Additionally, for any $\alpha \in \mathbb{R}$ and $j=0,\dots,k-1$,
  \begin{align*}
    g(x+\alpha p_{k-1}^N(x))^{T} q_j &=(g(x) -  \alpha\sum_{i=0}^{k-1}
    \frac{g(x)^{T} q_i}{q_i^{T}Hq_i} Hq_i)^{T} q_j \\
                                   &= g(x)^{T} q_j -\alpha \frac{g(x)^{T}
                                     q_j}{q_j^{T}Hq_j}q_j^{T}Hq_j
    = (1-\alpha) g(x)^{T} q_j,
                                     \end{align*}
where the conjugacy of the vectors $q_0, \dots, q_{k-1}$ has been used.
\end{proof}

For a point $x \in x_0 + \K_k(g_0,H)$, the subspace Newton step to $\hat x_{k}$ may therefore be recovered using a Newton scaled conjugate basis for $\K_k(g_0,H)$. Once the step $p_{k-1}^{N}(x)$ is known explicitly, it is possible to recover $\hat x_k$ and therefore the gradient $\hat g_k$. Since $\hat g_k$ is orthogonal to $\K_{k}(g_0,H)$, the vectors $\{\hat g_k, q_{k-1}, \dots, q_0\}$ form a basis for $\K_{k+1}(g_0,H)$. Using such a basis, we may construct the subspace Newton step $p_{k}^{N}(x)$ to the minimizer $\hat x_{k+1} \in x_0 + \K_{k+1}(g_0,H)$ for any $x \in x_0 + \K_k(g_0,H)$.

\begin{proposition}\label{sub:newton}
For $k < r$, let $x \in x_0 + \K_{k}(g_0,H)$ and let $p_{k-1}^{N}(x)$ be the subspace Newton step to $\hat x_k$. If $q_{k-1} \parallel \hat q_{k-1}$ and $\hat g_k$ is the gradient at $\hat x_k$, then the subspace Newton step $p_{k}^{N}(x)$  to $\hat x_{k+1}$ is
\begin{align*}p_{k}^{N}(x) & = \hat q_k+ p_{k-1}^N(x)
  \end{align*}
 with
 \[\hat q_k =  \frac{1}{\widehat{\sigma}_{k}}\left(-\hat g_k + \frac{\hat g_k^{T} H q_{k-1}}{q_{k-1}^{T}Hq_{k-1}}q_{k-1}\right) \text{and } \widehat{\sigma}_{k} = \frac{ \hat g_k^{T} H \hat g_k q_{k-1}^{T} H q_{k-1} - (\hat g_k^{T} H q_{k-1})^2 }{\hat g_k^{T} \hat g_k  q_{k-1}^{T} H q_{k-1}}.\]
\end{proposition}

\begin{proof}
Since $\hat x_{k+1} \in x_0 + \K_{k+1}(g_0,H)$, the subspace Newton step $p_{k}^{N}(x)$ may be written as a linear combination of $ \hat g_k, q_{k-1}, \dots, q_0$ where $q_i \parallel \hat q_i$ for $i = 0, \dots, k-1$ which forms a basis for $\K_{k+1}(g_0,H)$.
Taking into account that $\hat g_k^{T} H q_i=0$, $i\le k-2$, and $q_i^{T} H
q_j=0$, $i\ne j$, yields
\[
  \mtx{ccccc}{ \hat g_k^{T} H \hat g_k & \hat g_k^{T} H q_{k-1} & 0 & 0 & 0 \\
                      q_{k-1}^{T} H \hat g_k & q_{k-1}^{T} H q_{k-1} & 0 &
                      0 & 0 \\
                      0 & 0 & q_{k-2}^{T} H q_{k-2} & 0  & 0 \\
                      0 & 0 & 0 & \ddots & 0 \\
                      0 & 0 & 0 & 0 & q_0^{T} H q_0 }
  \mtx{c}{\beta_k \\ \beta_{k-1} \\ \beta_{k-2} \\ \vdots \\ \beta_0 }
  =
  -   \mtx{c}{ \hat g_k^{T} g(x) \\ q_{k-1}^{T} g(x) \\ q_{k-2}^{T} g(x) \\ \vdots \\
    q_0^{T} g(x) }.
\]
Then, \[\beta_i = -\frac{g(x)^{T}q_i}{q_i^{T}Hq_i}, \qquad i = 0, \dots, k-2\] is the Newton scaling with respect to $x$.

Noting that $\hat g_k = g(x)+ H p_{k-1}^{N}(x)$, we may write
$g(x) = \hat g_k - H p_{k-1}^{N}(x)$. Since
\[p_{k-1}^N(x) = \sum_{i=0}^{k-1}-\frac{g(x)^{T}q_i}{q_i^{T}Hq_i}q_i.
  \] by Lemma~\ref{newton:step}, we have
\[
  -\hat g_k^{T} g(x)  =  -\hat g_k^{T} (\hat g_k - H p_{k-1}^{N}(x))=
  - \hat g_k^{T} \hat
  g_k - \frac{g(x)^{T}q_{k-1}}{q_{k-1}^{T}Hq_{k-1}}\hat g_k^{T} H q_{k-1},
\]
where the orthogonality of $\hat g_k$ to $H
q_i, i = 0, \dots, k-2$ is used.

Then, $\beta_k$ and $\beta_{k-1}$ are given by
\[
  \mtx{ccccc}{ \hat g_k^{T} H \hat g_k & \hat g_k^{T} H q_{k-1} \\
                      q_{k-1}^{T} H \hat g_k & q_{k-1}^{T} H q_{k-1} }
  \mtx{c}{\beta_k \\ \beta_{k-1} }
  =
  \mtx{c}{- \hat g_k^{T} \hat
  g_k - \frac{g(x)^{T}q_{k-1}}{q_{k-1}^{T}Hq_{k-1}}\hat g_k^{T} H q_{k-1} \\
   -q_{k-1}^{T}g(x) }.
\]
Let $\beta_{k-1}= - \frac{g(x)^{T}q_{k-1}}{q_{k-1}^{T}Hq_{k-1}}+\gamma_{k-1}$, then this is reduced to
\[
  \mtx{ccccc}{ \hat g_k^{T} H \hat g_k & \hat g_k^{T} H q_{k-1} \\
                      q_{k-1}^{T} H \hat g_k & q_{k-1}^{T} H q_{k-1} }
  \mtx{c}{\beta_k \\ \gamma_{k-1} }
  =
  \mtx{c}{ -\hat g_k^{T} \hat g_k \\
   0 },
\]
so that
\begin{eqnarray*}
  \mtx{c}{\beta_k \\ \gamma_{k-1} } & = &
  - \frac{\hat g_k^{T} \hat g_k}
  { \hat g_k^{T} H \hat g_k q_{k-1}^{T} H q_{k-1} - (\hat g_k^{T} H q_{k-1})^2 }
  \mtx{ccccc}{ q_{k-1}^{T} H q_{k-1} & -\hat g_k^{T} H q_{k-1} \\
                      -q_{k-1}^{T} H \hat g_k & \hat g_k^{T} H \hat g_k }
                                              \mtx{c}{ 1 \\ 0 } \\
  & = &
  - \frac{\hat g_k^{T} \hat g_k}
  { \hat g_k^{T} H \hat g_k q_{k-1}^{T} H q_{k-1} - (\hat g_k^{T} H q_{k-1})^2 }
  \mtx{ccccc}{ q_{k-1}^{T} H q_{k-1} \\
                      -q_{k-1}^{T} H \hat g_k }.
\end{eqnarray*}

The subspace Newton step to $\hat x_{k+1}$ is therefore formed by
\[ p_{k}^{N}(x) = \beta_{k} \hat g_k + \gamma_{k-1} q_{k-1} +  \sum_{i=0}^{k-1}-\frac{g(x)^{T}q_i}{q_i^{T}Hq_i}q_i=  \frac{1}{\widehat{\sigma}_{k}}\left(-\hat g_k + \frac{\hat g_k^{T} H q_{k-1}}{q_{k-1}^{T}Hq_{k-1}}q_{k-1}\right) + p_{k-1}^N(x)\]
where
\[ \widehat{\sigma}_k = -\frac{1}{\beta_k} = \frac{ \hat g_k^{T} H \hat g_k q_{k-1}^{T} H q_{k-1} - (\hat g_k^{T} H q_{k-1})^2 }{\hat g_k^{T} \hat g_k  q_{k-1}^{T} H q_{k-1}}.\]
\end{proof}

If we consider the case when $x = \hat x_{k}$, then $p_{k-1}^{N}(\hat x_k) = 0$.
Therefore
\[ p_{k}^{N} (\hat x_k) =\frac{1}{\widehat{\sigma}_k} \left(-\hat g_k + \frac{\hat g_k^{T}Hq_{k-1}}{q_{k-1}^{T}Hq_{k-1}} q_{k-1}\right) = \hat q_k,\]
is parallel to the search direction obtained using CG, BFGS, and
memoryless BFGS and is Newton-scaled with respect to $\hat x_k$.

Given that information is known on how $H$ acts on the whole space $\K_{k+1}(g_0,H)$, it is possible to calculate the subspace Newton step to $\hat x_{k+1}$ using the vectors $\hat g_k, q_{k-1}, \dots, q_0$. Building on this framework, we present a subspace quasi-Newton step under the assumption that we have incomplete information on the action of $H$ on the Krylov subspace $\K_{k+1}(g_0,H)$.

\section{Subspace quasi-Newton methods}\label{sec:qnsubspace}

To recover the subspace Newton step from $x \in x_0 + \K_k(g_0,H)$ to $\hat x_{k+1}$ as presented in the previous section, we require knowledge of how $H$ acts on the full subspace $\K_{k+1}(g_0,H)$.
In particular, sufficient information is needed to calculate the Newton scaling of the conjugate vectors $q_0, \dots, q_{k-1}$ with respect to $x$ as well as the scalars $\beta_k, \gamma_{k-1}$.
Suppose instead that we have a quasi-Newton approximation $B_k \succ 0$ such that $B_k q_i = H q_i$ for $i = 0, \dots, k-1$.
Then the only term that is unknown is $\hat g_k^{T}H\hat g_k$, which must be replaced by the approximation $\hat g_k^{T} B_k \hat g_k$.

If such an approximation is used, it is possible to replace $H$ with $B_k$ and use the construction in Proposition~\ref{sub:newton} to generate a step $p_k(x) =q_k + p_{k-1}^{N}(x)$ where
\[q_k = \frac{\widehat{\sigma}_k}{\sigma_k} \hat q_k \text{for }
  \sigma_k=\frac{ \hat g_k^{T} B_k \hat g_k q_{k-1}^{T} H q_{k-1} - (\hat g_k^{T} H q_{k-1})^2 }{\hat g_k^{T}\hat g_k q_{k-1}^{T}Hq_{k-1}},
\]
taking into account that $B_k q_{i} = H q_{i}, i = 0, \dots, k-1$. However, this construction requires explicit knowledge of the gradient $\hat g_k$. Additionally, $\sigma_k$ is undefined when $k = r$.
To avoid this, we will build an approximation $B_k$ which acts on the Krylov subspace $\K_k(g_0,H)$ for any $k \leq r$ in such a way that for any $x \in x_0 + \K_k(g_0,H)$, solving the system of equations $B_k p_k(x) = - g(x)$ yields a step $p_k(x) = q_k + p_{k-1}^N(x)$ where $q_k \parallel \hat q_k$, without requiring knowledge of $\hat g_k$.

\begin{lemma}\label{full:rank}
For $k \leq r$, let $x \in x_0 + \K_{k}(g_0,H)$ and let $Q_{k-1}=( q_{k-1} \ q_{k-2} \ \dots \ q_0)$ for $q_i \parallel \hat q_i, i = 0, \dots, k-1$.
Assume that $B_k$ is given by
\begin{equation}\label{full:memory}
  B_k = \sigma_k (I - Q_{k-1} (Q_{k-1}^{T} Q_{k-1})\inv Q_{k-1}^{T}) +
  H Q_{k-1} (Q_{k-1}^{T} H Q_{k-1})\inv Q_{k-1}^{T} H,
\end{equation}
where $\sigma_k>0$.
Then $B_kp_k(x) = -g(x)$ has a unique solution given by
\begin{eqnarray}\label{full:step}
p_k(x) & =
q_k + p_{k-1}^{N}(x)
\end{eqnarray}
where, for $k < r$, \[q_k= \frac{1}{\sigma_k} \left(-\hat g_k + \frac{\hat g_k^{T} H q_{k-1}}
        {q_{k-1}^{T} H q_{k-1}} q_{k-1} \right) = \frac{\widehat{\sigma}_k}{\sigma_k} \hat q_k,\]
for $\widehat{\sigma}_k$ and $\hat q_k$ given by Proposition~\ref{sub:newton}.
As a consequence, when $\sigma_k = \widehat{\sigma}_k$, the solution is the subspace Newton step $p_{k}^{N}(x)$ to $\hat x_{k+1}$.
For $k =r$, $q_r = \hat q_r = 0$ and the solution $p_r(x) = p_{r-1}^{N}(x)$ is independent of the choice of $\sigma_r$. \end{lemma}

\begin{proof}
Since $\hat g_k$ is the gradient at $\hat x_k$, it must be orthogonal to the space spanned by the columns of $Q_{k-1}$.
Additionally, by construction the term
\[ q_k = \frac{1}{\sigma_k}\left(-\hat g_k + \frac{\hat g_k^{T} H q_{k-1}}
        {q_{k-1}^{T} H q_{k-1}} q_{k-1}\right) \]
is conjugate to the vectors $q_0, \dots, q_{k-1}$ with respect to $H$ and belongs to $\K_{k+1}(g_0,H)$.
Therefore $B_k q_k = -\hat g_k$.

By Lemma~\ref{newton:step}, $p_{k-1}^{N}(x)$ is a linear combination of the vectors $q_0, \dots, q_{k-1}$. Therefore $B_k p_{k-1}^{N}(x) = H p_{k-1}^{N}(x)$. Since $p_{k-1}^{N}(x)$ is the subspace Newton step to $\hat x_k$, the gradient at $x$ may be rewritten as $g(x) = \hat g_k - Hp_{k-1}^{N}(x).$

Applying $B_k$ to \eqref{full:step} therefore yields the result

\[
B_k p_k(x)
         = - \hat g_k +  Hp_{k-1}^{N}(x)= - g(x).
\]
As $Q_{k-1}$ has full column rank, and $\sigma_k > 0$, $B_k$ is strictly positive definite by Lemma~\ref{lem-Bposdef} and the solution $p_k(x)$ is therefore unique. For the case where $\sigma_k = \widehat{\sigma}_k$, the solution is exactly $p_{k}^{N}(x) = \hat q_k + p_{k-1}^{N}(x)$.

When $k = r$, the gradient $\hat g_r = 0$. Therefore, by definition, $q_r = \hat q_r = 0$. For any $\sigma_r > 0$,
\[ B_r p_{r-1}^{N}(x) = H p_{r-1}^{N}(x) = -g(x), \]
so $p_{r-1}^{N}(x) = p_r(x)$.
\end{proof}

The approximation $B_k$ given by \eqref{full:memory} requires a full basis of $\K_k(g_0,H)$ to calculate the step $p_k(x)$.
In the proof of Lemma~\ref{full:rank}, however, it is unnecessary for $B_k$ to act on each basis vector $q_0, \dots, q_{k-1}$ as $H$ individually.
The key components of the proof require only that $B_k$ acts as $H$ on $q_{k-1}$ and $p_{k-1}^{N}(x)$.

We now show that it suffices to replace the full basis of $\K_k(g_0,H)$ by the vectors $q_{k-1}$ and $p_{k-1}^{N}(x)$ in the construction of a Hessian approximation $B_k(x)$.
Then $B_k(x)$ need only act as $H$ on a subspace of dimension at most two, regardless of the size of $\K_k(g_0,H)$.

\begin{proposition}\label{rank:4}
For $k \leq r$, let $x \in x_0 + \K_{k}(g_0,H)$, let $p_{k-1}^{N}(x)$ be the subspace Newton step to $\hat x_k$ and let $q_{k-1} \parallel \hat q_{k-1}$. For notational simplicity, the dependence of $P_{k-1}$ on $x$ is suppressed. Define
\[
P_{k-1} =
\begin{cases}
\begin{bmatrix} p_{k-1}^{N}(x) & q_{k-1} \end{bmatrix}, & \text{if } \rank \left(\begin{bmatrix} p_{k-1}^{N}(x) & q_{k-1} \end{bmatrix}\right) = 2, \\[10pt]

\begin{bmatrix} q_{k-1} \end{bmatrix}, & \text{otherwise.}
\end{cases}
\]
Assume that $B_k(x)$ is given by
\begin{equation}\label{memory:two}
  B_k(x) = \sigma_k (I - P_{k-1}(P_{k-1}^{T} P_{k-1})\inv P_{k-1}^{T}) +
  H P_{k-1}(P_{k-1}^{T} H P_{k-1})\inv P_{k-1}^{T} H.
\end{equation}
Then $B_k(x)p_k(x) = -g(x)$ has a unique solution
\begin{eqnarray}\label{rank4:step}
p_k(x) & = q_k + p_{k-1}^{N}(x)
\end{eqnarray}
which coincides with the solution in Lemma~\ref{full:rank} calculated by the approximation $B_k$ given by \eqref{full:memory}.
\end{proposition}

\begin{proof}
Since $k \leq r$, $q_{k-1}$ must be nonzero. Therefore $P_{k-1}$ has full column rank by definition, so $B_k(x)$ is well defined. Additionally, $B_k(x)$ is positive definite by Lemma~\ref{lem-Bposdef}, so a solution to the equation $B_k(x)p_k(x) = -g(x)$ must exist and be unique.

Let
\[ q_k = \frac{1}{\sigma_k} \left( - \hat g_k + \frac{\hat g_k^{T}Hq_{k-1}}{q_{k-1}^{T}Hq_{k-1}}q_{k-1}\right).\]
By construction, $q_k$ is conjugate to $\K_k(g_0,H)$ with respect to $H$. By Lemma~\ref{newton:step} this implies that $q_k^{T}Hp_{k-1}^{N}(x) = 0$. Therefore $B_k(x)q_k=-\hat g_k$. Additionally, by design $B_k(x)p_{k-1}^{N}(x) = Hp_{k-1}^{N}(x)$. Applying $B_k(x)$ to \eqref{rank4:step} yields
\[B_k(x)p_k(x) = -\hat g_k + Hp_{k-1}^{N}(x) = -g(x) \]
which follows from $p_{k-1}^{N}(x)$ being the subspace Newton step to $\hat x_k$.

Therefore $p_k(x) = q_k + p_{k-1}^{N}(x)$ is the unique solution to the system of equations  $B_k(x)p_k(x) = -g(x)$ and coincides with the solution calculated by the approximation $B_k$ given by \eqref{full:memory}.
\end{proof}

For any $x \in x_0 + \K_k(g_0,H)$ we can now construct a Hessian approximation using the subspace Newton step $p_{k-1}^{N}(x)$ and a vector $q_{k-1}\parallel \hat q_{k-1}$. This Hessian approximation along with the gradient yields a quasi-Newton step that is a combination of $p_{k-1}^N(x)$ and $q_k \parallel \hat q_k$. We have shown that there is a unique value $\widehat{\sigma}_k$ such that the resulting step is exactly the subspace Newton step $p_{k}^{N}(x)$ to $\hat x_{k+1}$. However, even if the resulting step is not the subspace Newton step, the point $x + \alpha p_k(x)$ for any step size $\alpha$ belongs to the affine space $x_0 + \K_{k+1}(g_0,H)$. We may use this property to create an iterative method which successively generates quasi-Newton steps regardless of step size.

\section{An exact quasi-Newton method without exact line search}\label{sec:qn}

Given a point $x \in x_0 + \K_{k-1}(g_0,H)$, generating the subspace quasi-Newton step as described in Proposition~\ref{rank:4} requires explicit knowledge of the direction of the subspace Newton step $p_{k-1}^{N}(x)$ from $x$ to $\hat x_k$ as well as the direction $\hat q_{k-1}$.
In addition, constructing the approximation requires the action of $H$ on these vectors to build the Hessian approximation $B_k$ given by \eqref{memory:two}.
We propose an iterative construction that produces subspace quasi-Newton steps for arbitrary step sizes. After finitely many iterations, it produces the full Newton step; a unit step then reaches the minimizer.
In this Hessian-based formulation, $\alpha_k$ may be chosen arbitrarily, including zero. In the first-order formulation of Subsection~\ref{subsec:firstorder}, it must be nonzero.
This method is described in Algorithm~\ref{alg:main}. A hat distinguishes quantities determined by the Krylov-space minimizers from freely chosen algorithmic quantities; in particular, $\widehat{\sigma}_k$ is determined by the quadratic problem, whereas $\sigma_k>0$ is a free parameter. Within the algorithm, $p_{k-1}^{N}$ denotes the stored subspace Newton step $p_{k-1}^{N}(x_k)$; for $k\geq r$, it denotes the full Newton step $x^{*}-x_k$.

Theoretically, if we have explicit knowledge of $H$ then even a zero step will lead to generation of the Newton step because new conjugate directions can be generated iteratively and Newton-scaled at the fixed point.
Without explicit knowledge of $H$, a nonzero step is required to evaluate a new gradient.
It is then possible to use a difference of gradients to learn how $H$ acts on the step taken. This is described in Subsection~\ref{subsec:firstorder}.

\begin{algorithm}[H]
\caption{A two-vector quasi-Newton construction for \eqref{qp}}\label{alg:main}
\begin{algorithmic}
\Require $x_{0} \in \rn, p_{-1}^{N} = 0, \sigma_0 > 0$
\State $k \gets 0$; \quad $g_k \gets Hx_k + c$;
\While{$\norm{g_k}\ne 0$}
\If {$k=0$}
\State $B_k \gets \sigma_k I$;
\Else
    \State Choose some $\sigma_{k} > 0$;
    \State $B_{k} \gets \sigma_{k} ( I - P_{k-1}(P_{k-1}^{T}P_{k-1})^{-1}P_{k-1}^{T}) + HP_{k-1}(P_{k-1}^{T}HP_{k-1})^{-1}P_{k-1}^{T}H$;
    \EndIf
    \State Solve $B_kp_k = -g_{k}$;
    \State Choose step size $\alpha_k$;
    \State $x_{k+1} \gets x_{k} + \alpha_k p_{k}$;
    \State $g_{k+1} \gets g_{k} + \alpha_k Hp_{k}$;
    \If{$\norm{g_{k+1}} = 0$}
        \State \textbf{break}
    \EndIf
    \State $q_k \gets p_k - p_{k-1}^{N}$;
    \If{$\norm{q_k} \ne 0$}
        \State $p^{N}_k \gets \left(-\frac{g_k^{T}q_k}{q_k^{T}Hq_k} - 1\right)q_k + (1-\alpha_k) p_k$;
        \If{$\rank( \begin{bmatrix} p_{k}^{N} & q_k\end{bmatrix} )=2$}
            \State $P_k \gets \begin{bmatrix} p_{k}^{N} & q_k\end{bmatrix}$;
        \Else
            \State $P_k \gets \begin{bmatrix} q_k \end{bmatrix}$;
        \EndIf
    \Else
        \State $p_{k}^{N} \gets (1-\alpha_k)p_{k-1}^{N}$;
        \State $P_k \gets \begin{bmatrix} p_{k}^{N} \end{bmatrix}$;
    \EndIf
    \State $k \gets k+1$;
\EndWhile
\end{algorithmic}
\end{algorithm}

\begin{theorem}\label{main:thm}
Let $x_0$ be fixed, and let $x_i$, $i=0,\dots,k$, be generated by
Algorithm~\ref{alg:main}, where the step lengths $\alpha_i$, $i<k$,
are arbitrary and $\sigma_i>0$, $i\leq k$. For $k\leq r$,
Algorithm~\ref{alg:main} generates exactly the step $p_k(x_k)$ defined
in Proposition~\ref{rank:4}. For $k\geq r$, the generated direction is
the Newton step at $x_k$. Consequently, once $k\ge r$, choosing
$\alpha_k=1$ terminates the algorithm.
\end{theorem}

\begin{proof}
When $k < r$, to show that the search direction $p_{k+1}$ generated by Algorithm~\ref{alg:main} is exactly $p_{k+1}(x_{k+1})$ given by \eqref{rank4:step}, it is sufficient to show that $q_{k} \parallel \hat q_{k}$ and that $p_{k}^N = p_{k}^N(x_{k+1})$. We proceed by induction. For $k=0$, clearly $\sigma_0q_0 = -g(x_0)$ and $p_0^N$ is Newton scaled with respect to $x_1$.

Assume that $q_{k-1}\parallel \hat q_{k-1}$ and $p_{k-1}^{N} = p_{k-1}^{N}(x_{k})$, then by Proposition~\ref{rank:4} $p_{k}$ is exactly $p_{k}(x_{k})$ which is of the form $p_{k} = q_{k} + p_{k-1}^{N}$, therefore $q_{k} \parallel \hat q_{k}$.

As the search directions generated by the algorithm coincide with those given in Proposition~\ref{rank:4} for $i \leq k$, it is straightforward to show that $x_{k+1} \in x_0 + \K_{k+1}(g_0,H)$. Using the formulation for the subspace Newton step in Lemma~\ref{newton:step}, we have
\begin{align}
  p_k^N(x_{k+1}) & =
                   -  \sum_{i=0}^k \frac {g(x_{k+1})^{T} q_i}{q_i^{T} H
                   q_i} q_i \nonumber\\
  & = - \sum_{i=0}^{k}\frac{(g_k+ \alpha_k Hq_k +\alpha_k H
    p_{k-1}^N(x_k))^{T}q_i}{q_i^{T} H q_i}q_i \nonumber\\\
  & = \left( - \frac{g_k^{T} q_k}{q_k^{T} H q_k } - \alpha_k \right) q_k -  \sum_{i=0}^{k-1} \frac {g(x_{k} + \alpha_kp_{k-1}^N(x_k))^{T} q_i}{q_i^{T} H
                   q_i} q_i \nonumber\\\
  & = \left( - \frac{g_k^{T} q_k}{q_k^{T} H q_k } - \alpha_k \right) q_k +
    p_{k-1}^N (x_{k} + \alpha_kp_{k-1}^N(x_k)) \nonumber\\\
  & = \left( -\frac{g_k^{T} q_k}{q_k^{T} H q_k } - \alpha_k \right) q_k +
    (1-\alpha_k) p_{k-1}^N (x_k) \nonumber\\\
  & = \left( - \frac{g_k^{T} q_k}{q_k^{T} H q_k } - 1 \right) q_k +
    (1-\alpha_k) p_k , \label{alt:basis}
\end{align}
taking the conjugacy of $q_k$ with $p_{k-1}^{N}(x_k), q_{k-1}, \dots, q_0$ into account, as well as the property of the subspace Newton step
\[ p_{k-1}^{N}(x_k + \alpha_kp_{k-1}^{N}(x_k)) = \hat x_k - x_k - \alpha_kp_{k-1}^{N}(x_k) = (1-\alpha_k)p_{k-1}^{N}(x_k).\]

Since $q_{k}$ is nonzero, $B_{k+1}$ is constructed by the algorithm exactly as $B_{k+1}(x_{k+1})$ is in Proposition~\ref{rank:4}. Therefore $p_{k+1}$ is the search direction $p_{k+1}(x_{k+1})$ given by \eqref{rank4:step}.

As a consequence of Proposition~\ref{rank:4}, $p_r = p_{r-1}^{N}(x_r)$ is the Newton step and $q_r = \hat q_r = 0$.
Therefore for any step size $\alpha_r$, $p_{r}^{N} = (1 - \alpha_r)p_{r-1}^{N}$ is the Newton step from $x_{r+1}$.
A second induction then shows that, for any $k > r$, if $p_{k-1}^{N}$ is the Newton step and $q_k = 0$, then $B_k$ defined by the algorithm will generate the Newton step. The algorithm will therefore terminate if $\alpha_k = 1$.
\end{proof}

As stated in Theorem~\ref{main:thm}, the search direction $p_k$ generated by Algorithm~\ref{alg:main} for $k \geq r$ is precisely the Newton step.
Therefore, if a unit step size is chosen, the algorithm will terminate.
It is however possible to terminate in exactly $r$ iterations if and only if the Hessian approximation is built using the unique value $\widehat{\sigma}_{r-1}$ stated in Proposition~\ref{sub:newton}.

\begin{corollary}\label{corr:unit}
For the unique value $\widehat{\sigma}_{r-1}$ stated in Proposition~\ref{sub:newton}, Algorithm~\ref{alg:main} will generate the Newton step $p_{r-1} = p_{r-1}^{N}(x_{r-1})$, and the algorithm will terminate at iteration $r$ if $\alpha_{r-1} = 1$.
For all other values of $\sigma_{r-1}$, the algorithm will terminate only when $\alpha_k = 1$ for some $k \geq r$.
\end{corollary}

\begin{proof}
As stated in Theorem~\ref{main:thm}, for $k \leq r$ the search direction $p_k$ coincides with the search direction $p_k(x_k)$ from Proposition~\ref{rank:4}. Therefore, for the unique value of $\widehat{\sigma}_{r-1}$ stated in Proposition~\ref{sub:newton}
\[ p_{r-1} = \frac{\widehat{\sigma}_{r-1}}{\widehat{\sigma}_{r-1}} \hat q_{r-1} + p_{r-2}^{N}(x_{r-1}) = p_{r-1}^{N}(x_{r-1}). \]
Hence the algorithm will terminate at iteration $r$ if $\alpha_{r-1} = 1$. For all other values of $\sigma_{r-1}$, the search direction $p_{r-1}$ cannot be parallel to $p_{r-1}^{N}(x_{r-1})$ as $\hat q_{r-1} \ne 0$. Therefore the algorithm cannot terminate. Applying Theorem~\ref{main:thm}, for $k \geq r$ the search direction $p_k = p_{r-1}^{N}(x_k)$ is the Newton step, so the algorithm will terminate only if $\alpha_k = 1$.
\end{proof}

At each iteration of Algorithm~\ref{alg:main}, the Hessian approximation $B_k$ is exactly that described by Proposition~\ref{rank:4} for $k\leq r$ since it is built using the two vectors $q_{k-1}$ and $p_{k-1}^{N}(x_k)$.
However, the use of those two specific vectors themselves is not necessary.
What matters is that $B_k$ agrees with $H$ on the subspace spanned by the vectors $q_{k-1}$ and $p_{k-1}^{N}(x_k)$.
The Hessian approximation may therefore be built using any vectors which span the same space.
From $\eqref{alt:basis}$ it is straightforward to show that $p_{k-1}^{N}(x_k)$ may be replaced with $p_{k-1}$ to form an alternate basis for the subspace in question.
Calculation of the subspace Newton step is still used at each iteration to recover the direction $q_k$.

This method will generate the Newton step in a finite number of iterations regardless of step sizes chosen.
The action of $B_k$ on the correct subspace of two vectors is the vital component which yields the subspace quasi-Newton step of Proposition~\ref{rank:4}. As shown in Lemma~\ref{lem-Bsolve} solving $B_kp_k = -g_k$ can be done by solving two systems of linear equations where the matrix is $P_k^{T}HP_k$, whose dimension is at most $2 \times 2$.

\subsection{A first-order formulation of Algorithm~\ref{alg:main}}\label{subsec:firstorder}

The construction of $B_{k+1}$ in Algorithm~\ref{alg:main} requires information on how the Hessian $H$ acts on the vectors $q_k$ and $p_{k}^{N}$.
However, explicit knowledge of $H$ is unnecessary.
The Hessian-based formulation of Algorithm~\ref{alg:main} allows arbitrary step
lengths, including zero. A first-order implementation, however,
requires $\alpha_k\ne0$ so that a new gradient can be
evaluated.
When a nonzero step is taken at each iteration, the required Hessian
action can be recovered from a gradient difference.

\begin{lemma}\label{first:order}
For $k < r$, if $\alpha_k \neq 0$ is chosen in Algorithm~\ref{alg:main}, then the action of $H$ on $q_{k}$ can be calculated by
\begin{equation}\label{hq}
Hq_k = \frac{1}{\alpha_k}(g_{k+1} - g_k) - Hp_{k-1}^{N},
\end{equation}
where $Hp_{k-1}^{N} = B_kp_{k-1}^{N}.$
Using this formulation, $p_{k}^{N}$ and $Hp_{k}^{N}$ can be recovered by
\begin{gather}
p_{k}^{N} = \left(-\frac{\alpha_k g_k^{T}q_k}{q_k^{T}(g_{k+1}-g_k)} - 1\right)q_k + (1-\alpha_k)p_k \label{p} \\
Hp_{k}^{N} = \left(-\frac{\alpha_k g_k^{T}q_k}{q_k^{T}(g_{k+1}-g_k)} - 1\right)Hq_k + \left(\frac{1}{\alpha_k}-1\right)(g_{k+1}-g_k). \label{hp}
\end{gather}
When $k \geq r$, $p_{k}^{N} = (1 - \alpha_k)  p_k$, so
\begin{equation}\label{hp2}
Hp_{k}^{N} = \left(\frac{1}{\alpha_k} - 1\right)(g_{k+1}-g_k).
\end{equation}
\end{lemma}

\begin{proof}
The step generated by Algorithm~\ref{alg:main} is of the form $p_k = q_k + p_{k-1}^{N}$. Assuming a nonzero step size $\alpha_k$ is used so that $x_{k+1} = x_k + \alpha_k p_k$, evaluating the gradient $g_{k+1}$ yields
\[ \frac{1}{\alpha_k}(g_{k+1} - g_k) = Hp_k = Hq_k + Hp_{k-1}^{N}.\]
Given that the matrix $B_k$ acts as $H$ on $p_{k-1}^{N}$, the formulation of $Hq_k$ in \eqref{hq} follows directly.

As $p_{k-1}^{N} \in \K_k(g_0,H)$, it is conjugate to $q_k$ with respect to $H$. Therefore we may replace the term $q_k^{T}Hq_k$ by
\begin{align*}
q_k^{T}Hq_k &= q_k^{T}\left(\frac{1}{\alpha_k}(g_{k+1} - g_k) - Hp_{k-1}^{N}\right) \\
&= \frac{1}{\alpha_k}q_k^{T}(g_{k+1} - g_k)
\end{align*}
in \eqref{alt:basis} to obtain the expression \eqref{p}. The formulations for $Hp_{k}^{N}$ in \eqref{hp} and \eqref{hp2} follow immediately.
\end{proof}

Using the expressions for $p_{k}^{N}, Hq_k$ and $Hp_{k}^{N}$ presented in the lemma above means that Algorithm~\ref{alg:main} can be implemented using first-order information only to recover the action of $H$ on the necessary vectors.
The termination guarantees of the algorithm remain unchanged with this formulation.

\section{Conclusion}\label{sec:conclusion}

We have investigated the curvature information sufficient for finite Newton-step generation by quasi-Newton methods on strictly convex quadratic problems. In the classical exact-line-search setting, memoryless BFGS achieves finite termination using curvature information associated with one direction. Our results show that, when the preceding step lengths are arbitrary, information associated with at most one additional direction is sufficient.

More precisely, it is unnecessary for the quasi-Newton approximation to agree with the Hessian on the full Krylov subspace. Agreement on a suitably chosen subspace of dimension at most two suffices to generate a step comprising a subspace Newton step and a component parallel to the next conjugate direction. Applied recursively over the expanding Krylov subspaces, this construction generates the full Newton step after finitely many iterations; a subsequent unit step reaches the minimizer. The required Hessian actions can also be recovered from gradient differences, giving a first-order formulation.

The construction serves to establish this result on sufficient
curvature information and is analyzed in exact arithmetic. We do not propose it as a practical alternative to the
method of conjugate gradients. The main conclusion is instead structural: removing exact
line search requires only curvature information associated with one
additional direction, rather than agreement with the Hessian on the
full accumulated Krylov subspace.

\appendix
\section{Some linear algebra results}

This appendix collects two linear-algebra results. The first
lemma shows that the quasi-Newton matrices generated are positive
definite.

\begin{lemma}\label{lem-Bposdef}
  Let
  \[
    B=\sigma ( I - P(P^{T}P)^{-1}P^{T}) + HP(P^{T}HP)^{-1}P^{T} H,
  \]
  where $H$ is a symmetric positive definite $n\times n$ matrix, $P$
  is an $n\times m$ matrix of rank $m$, and $\sigma$ is a positive
  constant. Then, $B$ is positive definite.
\end{lemma}

\begin{proof}
  We have $B=B_1+B_2$, for
  \[
  B_1=\sigma ( I - P(P^{T}P)^{-1}P^{T}) \text{and} B_2=HP(P^{T}HP)^{-1}P^{T}
  H.
  \]
  Both $B_1$ and $B_2$ are symmetric and positive semidefinite, as
  $\sigma>0$ and $H$ is symmetric and positive definite. Their sum is
  positive definite if there is no nontrivial intersection of their
  nullspaces. We have $\Null(B_1)=\Range(P)$ in addition to
  $B_2P=HP$. As $H$ is positive definite and $P$ has full column rank,
  the nullspaces of $B_1$ and $B_2$ have a trivial intersection.
\end{proof}

The following lemma reduces a system with the low-rank
quasi-Newton matrix.

\begin{lemma}\label{lem-Bsolve}
  Let
  \[
    B=\sigma ( I - P(P^{T}P)^{-1}P^{T}) + HP(P^{T}HP)^{-1}P^{T} H,
  \]
  where $H$ is a symmetric positive definite $n\times n$ matrix, $P$
  is an $n\times m$ matrix of rank $m$, and $\sigma$ is a positive
  constant. Then, $Bp=-g$ if and only if
\begin{eqnarray*}
  P^{T} H P \beta & = & - P^{T} g, \\
P^{T} H P \delta & = & -P^{T} H g - ( \sigma P^{T} H P + P^{T} H^2 P)\beta, \\
  p & = & -\frac1{\sigma} \left(g + P \delta + HP \beta\right).
\end{eqnarray*}
\end{lemma}

\begin{proof}
  Lemma~\ref{lem-Bposdef} shows that $B$ is positive definite, so that
  $p$ is well defined and unique.  In addition, we have
\[
    (\sigma ( I - P(P^{T}P)^{-1}P^{T}) + HP(P^{T}HP)^{-1}P^{T} H ) p = -g
\]
if and only if
\[
  \mtx{ccc}{ \sigma I & P & HP \\
    P^{T} & \frac1\sigma P^{T} P & 0 \\
    P^{T} H & 0 & -P^{T} H P }
  \mtx{c}{ p \\ \delta \\ \beta }
  = -\mtx{c}{ g \\ 0 \\ 0 }
\]
for some $m$-dimensional vectors $\delta$ and $\beta$. Elimination of
$p$ in the second and third block of equations gives
\[
  \mtx{ccc}{ \sigma I & P & HP \\
    0 & 0 & -\frac1\sigma P^{T} H P  \\
    0  &  -\frac 1{\sigma} P^{T} H P  & -P^{T} H P - \frac1{\sigma} P^{T} H^2 P }
  \mtx{c}{ p \\ \delta \\ \beta }
  = \mtx{c}{ -g \\ \frac1\sigma P^{T} g \\ \frac1{\sigma} P^{T} H g }.
\]
We may therefore solve in turn
\begin{eqnarray*}
  P^{T} H P \beta & = & - P^{T} g, \\
P^{T} H P \delta & = & -P^{T} H g - ( \sigma P^{T} H P + P^{T} H^2 P)\beta, \\
  p & = & -\frac1{\sigma} \left(g + P \delta + HP \beta\right),
\end{eqnarray*}
completing the proof.
\end{proof}

Thus, a system with $B$ can be solved by solving two systems whose
order equals the number of columns of $P$. In the present
construction, these systems are of order at most two.

\section*{Use of generative AI}

The authors have used Microsoft Copilot during the preparation of this
manuscript to discuss alternative formulations and improve the clarity
and presentation of parts of the text. The authors have reviewed and
revised all suggested content and take full responsibility for the
final manuscript.

\printbibliography[heading=bibintoc, title={References}]

\end{document}